\documentclass[12pt]{article}
\usepackage{amssymb,amsmath}
\usepackage{hyperref}
\usepackage{graphicx}
\usepackage{color}
\usepackage{chngcntr}

\begin{document}

\title{\LARGE\bf The two-term Machin-like formula for pi with small arguments of the arctangent function}

\author{
\normalsize\bf S. M. Abrarov\footnote{\scriptsize{Dept. Earth and Space Science and Engineering, York University, Toronto, Canada, M3J 1P3.}}\, and B. M. Quine$^{*}$\footnote{\scriptsize{Dept. Physics and Astronomy, York University, Toronto, Canada, M3J 1P3.}}}

\date{April 15, 2017}
\maketitle

\begin{abstract}
In this paper we propose a new method for determination of the two-term Machin-like formula for pi with arbitrarily small arguments of the arctangent function. This approach excludes irrational numbers in computation and leads to a significant improvement in convergence with decreasing arguments of the arctangent function.\vspace{0.25cm}
\vspace{0.25cm}
\end{abstract}

\section{Methodology description}

In our recent publication we have derived a simple formula for pi \cite{Abrarov2016a} (see also \cite{OEIS} for the Mathematica code)
\begin{equation}\label{eq_1}
\pi = 2^{k+1}\arctan \left( \frac{\sqrt{2-a_{k-1}}}{a_k} \right),
\end{equation}
where the corresponding nested radicals are defined by recurrence relations ${{a}_{k}}=\sqrt{2+{{a}_{k-1}}}$,	${{a}_{1}}=\sqrt{2}$.
It is convenient to represent the equation \eqref{eq_1} as
\begin{equation}\label{eq_2}
\frac{\pi }{4}={{2}^{k-1}}\arctan \left( \frac{\sqrt{2-{{a}_{k-1}}}}{{{a}_{k}}} \right).
\end{equation}

The argument of the arctangent function in the formula \eqref{eq_2} is irrational as it is based on the nested radicals consisting of square roots of twos. As a consequence, the application of this formula requires multiple algebraic manipulations over a surd argument ${\sqrt{2-a_{k-1}}}/{a_k}$ of the arctangent function that causes some complexities in computing pi \cite{Guillera2016}. However, this problem may be effectively resolved by observing that the formula \eqref{eq_2} for pi can also be rewritten in form
\begin{equation}\label{eq_3}
\frac{\pi }{4}={{2}^{k-1}}\arctan \left( \frac{1}{{{u}_{1}}-\varepsilon } \right),
\end{equation}
where ${{u}_{1}}$ is assumed to be a positive rational number such that
\begin{equation}\label{eq_4}
{{u}_{1}}=\frac{{{a}_{k}}}{\sqrt{2-{{a}_{k-1}}}}+\varepsilon,	\qquad\qquad \left| \varepsilon  \right|<<{{u}_{1}}.
\end{equation}
Since the error term $\left| \varepsilon  \right|$ is significantly smaller than ${{u}_{1}}$, from equation \eqref{eq_3} it follows that
\[
\frac{\pi }{4}\approx {{2}^{k-1}}\arctan \left( \frac{1}{{{u}_{1}}} \right).
\]
Consequently, we may introduce a small remainder $\Delta $ in order to rearrange this equation into exact form
\begin{equation}\label{eq_5}
\frac{\pi }{4}={{2}^{k-1}}\arctan \left( \frac{1}{{{u}_{1}}} \right)+\Delta.
\end{equation}
Assuming now that
$$
\tan\left( \Delta  \right)=\frac{1}{{{u}_{2}}}\Leftrightarrow \Delta =\arctan \left( \frac{1}{{{u}_{2}}} \right)
$$
the equation \eqref{eq_5} can be expressed as
\begin{equation}\label{eq_6}
\frac{\pi }{4}={{2}^{k-1}}\arctan \left( \frac{1}{{{u}_{1}}} \right)+\arctan \left( \frac{1}{{{u}_{2}}} \right).
\end{equation}

The Machin-like formula for pi is given by \cite{Lehmer1938, Borwein1987}
$$
\frac{\pi }{4}=\sum\limits_{k=1}^{K}{{{\alpha }_{k}}\arctan \left( \frac{1}{{{\beta }_{k}}} \right)},
$$
where ${{\alpha }_{k}}$ and ${{\beta }_{k}}$ are either integers or, more generally, rationals. Consequently, the two-term Machin-like formula for pi can be represented in form
\begin{equation}\label{eq_7}
\frac{\pi }{4}={{\alpha }_{1}}\arctan \left( \frac{1}{{{\beta }_{1}}} \right)+{{\alpha }_{2}}\arctan \left( \frac{1}{{{\beta }_{2}}} \right).
\end{equation}

Consider the following identity \cite{Borwein1987}
$$
\arctan \left( x \right)=\frac{1}{2i}\ln \left( \frac{1+ix}{1-ix} \right)\Leftrightarrow \arctan \left( \frac{1}{x} \right)=\frac{1}{2i}\ln \left( \frac{x+i}{x-i} \right).
$$
Substituting this identity into the two-term Marchin-like formula \eqref{eq_7} for pi gives
$$
\frac{\pi }{4}=\frac{{{\alpha }_{1}}}{2i}\ln \left( \frac{{{\beta }_{1}}+i}{{{\beta }_{1}}-i} \right)+\frac{{{\alpha }_{2}}}{2i}\ln \left( \frac{{{\beta }_{2}}+i}{{{\beta }_{2}}-i} \right)
$$
or
\small
\[
\begin{aligned}
\frac{\pi }{2}i&=\ln \left( {{\left( \frac{{{\beta }_{1}}+i}{{{\beta }_{1}}-i} \right)}^{{{\alpha }_{1}}}} \right)+\ln \left( {{\left( \frac{{{\beta }_{2}}+i}{{{\beta }_{2}}-i} \right)}^{{{\alpha }_{2}}}} \right) \\
&=\ln \left( {{\left( \frac{{{\beta }_{1}}+i}{{{\beta }_{1}}-i} \right)}^{{{\alpha }_{1}}}}{{\left( \frac{{{\beta }_{2}}+i}{{{\beta }_{2}}-i} \right)}^{{{\alpha }_{2}}}} \right).
\end{aligned}
\]
\normalsize
Considering that ${{e}^{i\pi /2}}=i$ the exponentiation of this identity leads to the well-known equation \cite{Weisstein}
\begin{equation}\label{eq_8}
{{\left( \frac{{{\beta }_{1}}+i}{{{\beta }_{1}}-i} \right)}^{{{\alpha }_{1}}}}{{\left( \frac{{{\beta }_{2}}+i}{{{\beta }_{2}}-i} \right)}^{{{\alpha }_{2}}}}=i.
\end{equation}
Comparing equations \eqref{eq_6} with \eqref{eq_7}, from equation \eqref{eq_8} it follows now that
\begin{equation}\label{eq_9}
{{\left( \frac{{{u}_{1}}+i}{{{u}_{1}}-i} \right)}^{{{2}^{k-1}}}}\frac{{{u}_{2}}+i}{{{u}_{2}}-i}=i.
\end{equation}
Since in this equation $u_1$ is a rational number and $2^{k-1}$ is a positive integer, the real and imaginary parts
$$
{{\left( \frac{{{u}_{1}}+i}{{{u}_{1}}-i} \right)}^{{{2}^{k-1}}}}={{\left( \frac{u_{1}^{2}-1}{1+u_{1}^{2}}+i\frac{2{{u}_{1}}}{1+u_{1}^{2}} \right)}^{{{2}^{k-1}}}}
$$
are both rationals. This signifies that the real and imaginary parts 
$$
\frac{{{u}_{2}}+i}{{{u}_{2}}-i}=\frac{u_{2}^{2}-1}{1+u_{2}^{2}}+i\frac{2{{u}_{2}}}{1+u_{2}^{2}}
$$
must also be rationals. This is possible if and only if $u_2$ is a rational number.

It is not difficult to see that a solution of equation \eqref{eq_9} with respect to the unknown rational value ${{u}_{2}}$ is
\begin{equation}\label{eq_10}
{{u}_{2}}=\frac{2}{{{\left[ \left( {{u}_{1}}+i \right)/\left( {{u}_{1}}-i \right) \right]}^{{{2}^{k-1}}}}-i}-i.
\end{equation}
Thus, using equations \eqref{eq_4} and \eqref{eq_10} we can readily find arguments $1/{{u}_{1}}$ and $1/{{u}_{2}}$ of the arctangent function in the two-term Machin-like formula \eqref{eq_6} for pi. It should be noted that due to no restrictions these arguments may be chosen arbitrarily small.

\section{Implementation}

\subsection{Arctangent function}

The Maclaurin series expansion of the arctangent function, also known historically as the Gregory's series \cite{Lehmer1938, Borwein1989}, is given by
\begin{equation}\label{eq_11}
\arctan \left( x \right)=x-\frac{{{x}^{3}}}{3}+\frac{{{x}^{5}}}{5}-\frac{{{x}^{7}}}{7}+\cdots =\sum\limits_{n=1}^{\infty }{\frac{{{\left( -1 \right)}^{n+1}}}{2n-1}{{x}^{2n-1}}}.
\end{equation}
Since this equation can be expressed as
$$
\arctan \left( x \right)=x+O\left( {{x}^{3}} \right),
$$
one can see that due to vanishing term $O\left( {{x}^{3}} \right)$ the accuracy of an arctangent function improves as its argument $x$ decreases. This strongly motivated us to look for the Machin-like formula for pi with small arguments of the arctangent function.

We have shown previously that the arctangent function can be represented as \cite{Abrarov2016b}
\[
\arctan \left( x \right)=i\underset{M\to \infty }{\mathop{\lim }}\,\sum\limits_{m=1}^{\left\lfloor \frac{M}{2}+1 \right\rfloor }{\frac{1}{2m-1}\left( \frac{1}{{{\left( 1+2i/x \right)}^{2m-1}}}-\frac{1}{{{\left( 1-2i/x \right)}^{2m-1}}} \right)}.
\]
The derivation of this equation described in the work \cite{Abrarov2016b} is somehow tedious and based on the formula for numerical integration with enhanced midpoints in subintervals
\[
I = \int\limits_0^1 {f\left( t \right)dt}  = \underset{M\to \infty }{\mathop{\lim }}\sum\limits_{\ell  = 1}^L {\sum\limits_{m = 0}^M  {\frac{{{{\left( { - 1} \right)}^m} + 1}}{{{{\left( {2L} \right)}^{m + 1}}\left( {m + 1} \right)!}}{{\left. {{f^{\left( m \right)}}\left( t \right)} \right|}_{t = \frac{{\ell  - 1/2}}{L}}}} }
\]
where
\[
f\left( t \right) = \frac{x}{1+x^2 t^2}.
\]
However, Jes\'us Guillera has found recently an elegant and simple proof for this series expansion of the arctangent function  \cite{Guillera2016}.

\subsubsection*{Proof}

While $M$ tends to infinity, the upper bound $\left\lfloor M/2+1 \right\rfloor $ in summation also tends to infinity. Consequently, the series expansion of the arctangent function above can be simplified in form
\begin{equation}\label{eq_12}
\arctan \left( x \right)=i\sum\limits_{m=1}^{\infty }{\frac{1}{2m-1}\left( \frac{1}{{{\left( 1+2i/x \right)}^{2m-1}}}-\frac{1}{{{\left( 1-2i/x \right)}^{2m-1}}} \right)}.
\end{equation}

Consider the following identity \cite{Castellanos1988}
\begin{equation}\label{eq_13}
\arctan \left( c \right)+\arctan \left( d \right)=\arctan \left( \frac{c+d}{1-c\,d} \right).
\end{equation}
Assuming
$$
c=\frac{1}{2/x+i}
$$
and
$$
d=\frac{1}{2/x-i},
$$
it is easy to see that
$$
\frac{c+d}{1-c\,d}=x.
$$
Consequently, we can rearrange the identity \eqref{eq_13} in a reformulated form
\begin{equation}\label{eq_14}
\arctan \left( x \right)=\arctan \left( \frac{1}{2/x+i} \right)+\arctan \left( \frac{1}{2/x-i} \right).
\end{equation}
Applying now the Gregory's series \eqref{eq_11} with respect to the arguments $c$ and $d$
\footnotesize
\[
\arctan \left( c \right)=\sum\limits_{n=1}^{\infty }{\frac{{{\left( -1 \right)}^{n+1}}}{2n-1}{{c}^{2n-1}}\Leftrightarrow }\arctan \left( \frac{1}{2/x+i} \right)=\sum\limits_{n=1}^{\infty }{\frac{{{\left( -1 \right)}^{n+1}}}{2n-1}{{\left( \frac{1}{2/x+i} \right)}^{2n-1}}},
\]
\normalsize

\footnotesize
\[
\arctan \left( d \right)=\sum\limits_{n=1}^{\infty }{\frac{{{\left( -1 \right)}^{n+1}}}{2n-1}{{d}^{2n-1}}\Leftrightarrow \arctan \left( \frac{1}{2/x-i} \right)}=\sum\limits_{n=1}^{\infty }{\frac{{{\left( -1 \right)}^{n+1}}}{2n-1}{{\left( \frac{1}{2/x-i} \right)}^{2n-1}}}
\]
\normalsize
to both terms on the right side of equation \eqref{eq_14} immediately yields the series expansion of the arctangent function \eqref{eq_12}. This completes the proof.
\vspace{0.5cm}

The computational test we performed shows that the series expansion \eqref{eq_12} is more rapid in convergence than the Euler's formula \cite{Chien-Lih2005}
$$
\arctan \left( x \right)=\sum\limits_{n=0}^{\infty }{\frac{{{2}^{2n}}{{\left( n! \right)}^{2}}}{\left( 2n+1 \right)!}\frac{{{x}^{2n+1}}}{{{\left( 1+{{x}^{2}} \right)}^{n+1}}}}.
$$
In particular, with same number of the summation terms the series expansion \eqref{eq_12} is more accurate in computation by many orders of the magnitude than the Euler's formula and this tendency becomes especially evident when the argument $x$ tends to zero \cite{Guillera2016}.

Substituting series expansion \eqref{eq_12} into the two-term Machin-like formula \eqref{eq_6} for pi provides
\footnotesize
\[
\begin{aligned}
\frac{\pi}{4} = i&\sum\limits_{m=1}^{\infty }\frac{1}{2m-1} \times \\
&\left[ {{2}^{k-1}}\left( \frac{1}{{{\left( 1+2i\,{{u}_{1}} \right)}^{2m-1}}}-\frac{1}{{{\left( 1-2i\,{{u}_{1}} \right)}^{2m-1}}} \right)+\frac{1}{{{\left( 1+2i\,{{u}_{2}} \right)}^{2m-1}}}-\frac{1}{{{\left( 1-2i\,{{u}_{2}} \right)}^{2m-1}}} \right]
\end{aligned}
\]
\normalsize
or
\footnotesize
\begin{equation}\label{eq_15}
\hspace{-10.5cm}\pi \approx 4i\sum\limits_{m=1}^{m_{max}}\frac{1}{2m-1} \times
\end{equation}
\[
\hspace{0.9cm}\left[ {{2}^{k-1}}\left( \frac{1}{{{\left( 1+2i\,{{u}_{1}} \right)}^{2m-1}}}-\frac{1}{{{\left( 1-2i\,{{u}_{1}} \right)}^{2m-1}}} \right)+\frac{1}{{{\left( 1+2i\,{{u}_{2}} \right)}^{2m-1}}}-\frac{1}{{{\left( 1-2i\,{{u}_{2}} \right)}^{2m-1}}} \right]\hspace{-0.12cm},
\]
\normalsize
where $m_{max}>>1$ is the truncating integer. Further, we will use this equation in order to estimate the convergence rate in computing pi at given rational values ${{u}_{1}}$ and ${{u}_{2}}$.

\subsection{Computation}

Generally, it is very difficult to find just by guessing a combination of two simultaneously small rational arguments $1/{{\beta }_{1}}$ and $1/{{\beta }_{2}}$ of the arctangent function in the two-term Machin-like formula \eqref{eq_7} for pi. For example, substituting three random integers, say ${{\beta }_{1}}={{10}^{9}}$, ${{\alpha}_{1}}=7$ and ${{\alpha}_{2}}=1$, into equation \eqref{eq_8} we obtain the following solution for the unknown value
\footnotesize
\[
\begin{aligned}
{{\beta }_{2}}&=\frac{1000000006999999978999999965000000035000000020999999992999999999}{999999992999999979000000035000000034999999978999999993000000001} \\
&=1.00000001400000009800\ldots \,\, \left(\text{rational}\right).
\end{aligned}
\]
\normalsize
We can see now that this is not an optimal way for computation since only the first argument $1/{{\beta }_{1}}={{10}^{-9}}$ of the arctangent function is small whereas the second argument $1/{{\beta }_{2}}$ of the arctangent function is considerably larger and close to the unity. Therefore, due to relatively large value of the second argument $1/{{\beta }_{2}}$ of the arctangent function, we must not expect a rapid convergence in computation by substituting these values into the two-term Machin-like formula \eqref{eq_7} for pi.

In order to resolve this problem we can apply the proposed methodology based on the equations \eqref{eq_4}, \eqref{eq_6} and \eqref{eq_10}. We will consider some examples at $k$ equal to 2, 3, $5$, $10$, $17$ and $23$.

At $k=2$ the equation \eqref{eq_4} yields
$$
{{u}_{1}}=\frac{{{a}_{2}}}{\sqrt{2-{{a}_{1}}}}+\varepsilon =\frac{\sqrt{2+\sqrt{2}}}{\sqrt{2-\sqrt{2}}}+\varepsilon =2.41421356237309504880\ldots +\varepsilon
$$
If we take
$$
\varepsilon =-0.01421356237309504880\ldots \,\,\left( \text{irrational} \right),
$$
then we get
$$
{{u}_{1}}=2.4=\frac{24}{10}.
$$
Thus, substituting ${{u}_{1}}=24/10$ and $k=2$ into equation \eqref{eq_10} results in a solution ${{u}_{2}}=-239$. Consequently, the two-term Machin-like formula reads
\begin{equation}\label{eq_16}
\begin{aligned}
\frac{\pi }{4}& =2\arctan \left( \frac{10}{24} \right)+\arctan \left( \frac{1}{-239} \right) \\ 
 & =2\arctan \left( \frac{10}{24} \right)-\arctan \left( \frac{1}{239} \right).  
\end{aligned}
\end{equation}

At $k=3$ the equation \eqref{eq_4} yields
\[
\begin{aligned}
{{u}_{1}}& =\frac{{{a}_{3}}}{\sqrt{2-{{a}_{2}}}}+\varepsilon =\frac{\sqrt{2+\sqrt{2+\sqrt{2}}}}{\sqrt{2-\sqrt{2+\sqrt{2}}}}+\varepsilon \\
& =5.02733949212584810451\ldots +\varepsilon.
\end{aligned}
\]
Taking the error term as
$$
\varepsilon =-0.02733949212584810451\ldots  \,\,\left( \text{irrational} \right)
$$
leads to ${{u}_{1}}=5$. Substituting now ${{u}_{1}}=5$ and $k=3$ into equation \eqref{eq_10} we get again the same negative integer ${{u}_{2}}=-239$. Consequently, equation \eqref{eq_6} becomes
\begin{equation}\label{eq_17}
\begin{aligned}
\frac{\pi }{4}& = 4\arctan \left( \frac{1}{5} \right)+\arctan \left( \frac{1}{-239} \right) \\ 
& = 4\arctan \left( \frac{1}{5} \right)-\arctan \left( \frac{1}{239} \right).  
\end{aligned}
\end{equation}
This equation is well-known as the Machin\text{'}s formula for pi \cite{Lehmer1938, Borwein1987, Castellanos1988, Borwein2015, Borwein1989}.
 
It should be noted that comparing \eqref{eq_16} and \eqref{eq_17} we can find an interesting relation 
\[
4\arctan \left( \frac{1}{5} \right)=2\arctan \left( \frac{10}{24} \right).
\]

When $k=5$ the equation \eqref{eq_4} provides
\[
\begin{aligned}
{{u}_{1}}&=\frac{{{a}_{5}}}{\sqrt{2-{{a}_{5}}}}+\varepsilon =\frac{\sqrt{2+\sqrt{2+\sqrt{2+\sqrt{2+\sqrt{2}}}}}}{\sqrt{2-\sqrt{2+\sqrt{2+\sqrt{2+\sqrt{2}}}}}}+\varepsilon \\
&=20.35546762498718817831\ldots +\varepsilon.
\end{aligned}
\]
By taking
$$
\varepsilon =-0.35546762498718817831\ldots \,\,\left( \text{irrational} \right)
$$
we have ${{u}_{1}}=20$. Consequently, substituting ${{u}_{1}}=20$ and $k=5$ into equation \eqref{eq_10} we obtain
\[
\begin{aligned}
{{u}_{2}}=&-\frac{945426570789006031681}{13176476709447727679} \\
=&-71.75109034353024503462 \ldots \,\, \left(\text{rational}\right).
\end{aligned}
\]

Applying $k=10$ in equation \eqref{eq_4} leads to
$$
{{u}_{1}}=\frac{{{a}_{10}}}{\sqrt{2-{{a}_{9}}}}+\varepsilon =651.89813557739378661810\ldots +\varepsilon.
$$
Taking the error term as
$$
\varepsilon =-0.89813557739378661810\ldots \,\,\left( \text{irrational} \right)
$$
one obtains ${{u}_{1}}=651$. With ${{u}_{1}}=651$ and $k=10$ the equation \eqref{eq_10} yields
\[
\begin{aligned}
{{u}_{2}} & =-\frac{\overbrace{4370834256\ldots 5125120001}^{1,364\,\,\text{digits}}}{\underbrace{4736031894\ldots 5364787199}_{1,361\,\,\text{digits}}} \\
& =-922.88953146392823766085\ldots \,\,\left( \text{rational} \right).
\end{aligned}
\]

For $k=17$ the equation \eqref{eq_4} provides
$$
{{u}_{1}}=\frac{{{a}_{17}}}{\sqrt{2-{{a}_{16}}}}+\varepsilon =83443.02679976888016443942\ldots +\varepsilon.
$$
With error term taken as
$$
\varepsilon =0.02679976888016443942\ldots \,\,\left( \text{irrational} \right)
$$
we obtain ${{u}_{1}}=83443$. Consequently, substituting ${{u}_{1}}=83443$ and $k=17$ into equation \eqref{eq_10} gives
\[
\begin{aligned}
{{u}_{2}}& =-\frac{\overbrace{1617128975\ldots 8890856449}^{312,665\,\,\text{digits}}}{\underbrace{4079206389\ldots 1607609343}_{312,658\,\,\text{digits}}} \\
& =-3.96432252145804935647\ldots \times {{10}^{6}}\,\,\left( \text{rational} \right).
\end{aligned}
\]

Lastly, at $k=23$ equation \eqref{eq_4} results in
$$
{{u}_{1}}=\frac{{{a}_{23}}}{\sqrt{2-{{a}_{22}}}}+\varepsilon =5340353.71544080937733612922\ldots +\varepsilon
$$
and with
$$
\varepsilon =-0.01544080937733612922\ldots \,\,\left( \text{irrational} \right)
$$
we have
$$
{{u}_{1}}=5340353.7=\frac{53403537}{10}
$$
Substituting ${{u}_{1}}=53403537/10$ and $k=23$ into equation \eqref{eq_10} leads to
\[
\begin{aligned}
{{u}_{2}} & =-\frac{\overbrace{1009275657\ldots 1092218881}^{32,411,779\,\,\text{digits}}}{\underbrace{2291921400\ldots 3550735359}_{32,411,770\,\,\text{digits}}}\\
& =-4.40362247052490238495\ldots \times {{10}^{8}}\,\,\left( \text{rational} \right).
\end{aligned}
\]

The computational test reveals that with $k$ equal to $2$, $3$, $5$, $10$, $17$ and $23$ the equation \eqref{eq_15} contributes for $1$, $2$, $3$, $6$, $10$ and $14$ digits of pi, respectively, at each increment of the truncating integer $m_{max}$ by one. As we can see the convergence rate significantly improves with increasing integer $k$. In particular, at $k=23$ the convergence rate ($14$ digits per term) is approximately same as that of provided by the Chudnovsky formula for pi \cite{Baruah2009}.

Although equation \eqref{eq_10} is simple, due to rapid growth of the power ${{2}^{k-1}}$ the larger values of the integer $k$ require extensive computation. As a result, by using a typical desktop computer we were able to determine the rational numbers ${{u}_{2}}$ for the integers only up to $k=23$. However, there are no any theoretical limitations for the integer $k$ and we can estimate the convergence rate even at $k>23$ by using the identity \eqref{eq_2}. In particular, substituting equation \eqref{eq_12} into identity \eqref{eq_2} leads to
\small
\[
\begin{aligned}
\frac{\pi}{4}=&\,2^{k-1}i\sum\limits_{m=1}^{\infty }\frac{1}{2m-1} \times \\
&\left( \frac{1}{{{\left( 1+2i\,a_k/\sqrt{2-a_{k-1}} \right)}^{2m-1}}}-\frac{1}{{{\left( 1-2i\,a_k/\sqrt{2-a_{k-1}} \right)}^{2m-1}}} \right)
\end{aligned}
\]
\normalsize
or
\small
\begin{equation}\label{eq_18}
\begin{aligned}
\frac{\pi}{4}\approx &\,2^{k-1}i\sum\limits_{m=1}^{m_{max} }\frac{1}{2m-1} \times \\
&\left( \frac{1}{{{\left( 1+2i\,a_k/\sqrt{2-a_{k-1}} \right)}^{2m-1}}}-\frac{1}{{{\left( 1-2i\,a_k/\sqrt{2-a_{k-1}} \right)}^{2m-1}}} \right).
\end{aligned}
\end{equation}
\normalsize
It is easy to verify by straightforward substitution that at $k=40$ the equation \eqref{eq_18} provides $24$ digits of pi per increment of the truncating integer $m_{max}$ just by one. As it follows from the equation \eqref{eq_4}
$$
\frac{1}{u_1} \approx {\frac{\sqrt{2-a_{k-1}}}{a_k}}.
$$
By choosing a sufficiently small $\left|\varepsilon\right| << 1$, from equation \eqref{eq_10} we can always get such a value of $u_2$ that satisfies (see examples above)
$$
\frac{1}{\left|u_2\right|} < {\frac{1}{u_1}}\, .
$$
Consequently, due to smallness of the arguments $1/u_1$ and $1/u_2$ (by absolute value) of the arctangent function we expect that at $k=40$ the equation \eqref{eq_15} can also contribute to $24$ additional digits of pi per term incremented.

Thus, at $k>23$ the convergence rate can be significantly increased further by using more powerful computers for determination of the constant $u_2$ such that $\left|{{u}_{2}}\right|>>1$. Consequently, using this approach we can obtain arbitrarily small arguments $1/{{u}_{1}}$ and $1/{{u}_{2}}$ (by absolute value) of the arctangent function in order to improve convergence rate in computation involving the two-term Machin-like formula \eqref{eq_6} for pi.

\section*{Acknowledgments}

This work is supported by National Research Council Canada, Thoth Technology Inc. and York University. The authors wish to thank Dr. Jes\'us Guillera (University of Zaragoza) for the proof to equation \eqref{eq_12}, constructive discussions and useful information.

\bigskip

\end{document}